\documentclass[a4paper,12pt]{article}

\usepackage{pstricks}
\usepackage{graphicx,psfrag}
\usepackage{amsmath}
\usepackage{amssymb}
\usepackage{amsthm}
\usepackage{color}

      \newcommand {\al}   {\alpha}          \newcommand {\bt}  {\beta}
                
      \newcommand {\del}  {\delta}          
              \newcommand {\ve}   {\varepsilon}

      \newcommand {\om}   {\omega}          \newcommand {\Om}  {\Omega}
      \newcommand {\pl}   {\partial}

      \newcommand {\RRR}  {{\mathbb R}}

     \newcommand {\beq}  {\begin{equation}}
      \newcommand {\eeq}  {\end{equation}}
      
      \newcommand{\cl}{\text{cl}\,}

      \newtheorem{theorem}{Theorem}

      \newtheorem{zam}{Remark}

\title{A note on Newton's problem of minimal resistance\\ for convex bodies}

\author{
Alexander Plakhov\\
Center for R\&{}D in Mathematics and Applications, Department of Mathematics,\\ University of Aveiro, 3810-193, Portugal\\
and Institute for Information Transmission Problems, Moscow, 127051, Russia
}

\date{}

\begin{document}
\maketitle

\begin{abstract}
We consider the following problem: minimize the functional $\int_\Omega f(\nabla u(x))\, dx$ in the class of concave functions $u: \Omega \to [0,M]$, where $\Omega \subset \mathbb{R}^2$ is a convex body and $M > 0$. If $f(x) = 1/(1 + |x|^2)$ and $\Om$ is a circle, the problem is called Newton's problem of least resistance. It is known that the problem admits at least one solution. We prove that if all points of $\partial\Omega$ are regular and ${(1+|x|)f(x)}/(|y|f(y)) \to +\infty$ as $(1+|x|)/|y| \to 0$ then a solution $u$ to the problem satisfies $u\rfloor_{\partial\Omega} = 0$. This result proves the conjecture stated in 1993 in the paper by Buttazzo and Kawohl \cite{BK} for Newton's problem.
\end{abstract}

\begin{quote}
\end{quote}

\begin{quote}
\end{quote}

\section{Introduction}

Consider the following simple mechanical model. A convex body moves forward in a homogeneous medium composed of point particles. The medium is extremely rare, so as mutual interaction of particles is neglected. There is no thermal motion of particles, that is, the particles are initially at rest. When colliding with the body, each particle is reflected elastically. As a result of collisions, there appears the drag force that acts on the body and slows down its motion.

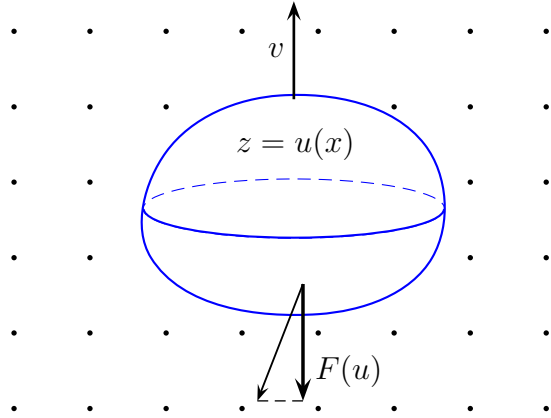
\begin{figure}[h]
\begin{picture}(0,165)
\scalebox{1}{
\rput(6.25,1.2){

\psellipse[linecolor=blue,linewidth=0.8pt](1.68,1.65)(1.99,0.4)
\pspolygon[fillstyle=solid,fillcolor=white,linewidth=0pt,linecolor=white](-0.4,1.65)(4,1.65)(4,2.5)(-0.4,2.5)
\psecurve[linecolor=blue,linewidth=0.8pt](-0.2,1)(2,0.25)(3.5,1)(3.4,2.5)(1.5,3.15)(0,2.4)(-0.2,1)(2,0.25)(3.5,1)
\psellipse[linecolor=blue,linewidth=0.4pt,linestyle=dashed](1.68,1.65)(1.99,0.4)
\psline[arrows=->,arrowscale=1.5,linewidth=1.3pt](1.8,0.65)(1.8,-0.9)
\psline[arrows=->,arrowscale=1.5,linewidth=0.7pt](1.8,0.65)(1.2,-0.9)
\psline[linestyle=dashed,linewidth=0.5pt](1.2,-0.9)(1.8,-0.9)
\rput(2.4,-0.5){$F(u)$}
\psline[arrows=->,arrowscale=1.5,linewidth=1pt](1.68,3.1)(1.68,4.4)
\rput(1.45,3.75){$v$}
\rput(1.7,2.5){$z = u(x)$}
\psdots[dotsize=2pt](-2,-1)(-2,0)(-2,1)(-2,2)(-2,3)(-2,4)
(-1,-1)(-1,0)(-1,1)(-1,2)(-1,3)(-1,4)
(0,-1)(0,0)(0,3)(0,4)
(1,-1)(1,0)(1,4)
(2,-1)(2,0)(2,4)
(3,-1)(3,0)(3,3)(3,4)
(4,-1)(4,0)(4,1)(4,2)(4,3)(4,4)
(5,-1)(5,0)(5,1)(5,2)(5,3)(5,4)

}}
\end{picture}
\caption{A convex body moving in a rarefied medium.}
\label{figRes}
\end{figure}

Take a coordinate system in $\RRR^3$ with the coordinates $x = (x^1, x^2), z$ connected with the body such that the $z$-axis is parallel and co-directional to the velocity of the body. Let the upper part of the body's surface be the graph of a concave function $u : \Om \to \RRR$, where $\Om$ is the projection of the body on the $x$-plane. Then the $z$-component of the drag force equals $-2\rho v^2 F(u)$, where $v$ is the scalar velocity of the body, $\rho$ is the density of the medium, and
\beq\label{resN}
F(u) = \int_\Om f(\nabla u(x))\, dx,
\eeq
where $f(x) = 1/(1 + |x|^2)$. $F(u)$ is called {\it resistance} of the body. See Fig.~\ref{figRes}.

\begin{zam}\label{z1}
In more realistic models, where the particle-body interaction is not perfectly elastic, the resistance also has the form \eqref{resN} with a different function $f$. Assume, for example, that the velocity of the reflected particle lies in the tangent plane to the body's surface at the point of impact and is co-directional with the projection of the velocity of incidence on this plane. Assume also that the moduli of the velocities of incidence and reflection coincide. Then we have $f(x) = 1 - |x|/\sqrt{1 + |x|^2}$.

Another generalization of the model concerns a medium with thermal motion of particles. The resistance in such a model is the sum of two integrals $\int_\Om f_1(\nabla u_1(x))\, dx + \int_\Om f_2(\nabla u_2(x))\, dx$, where the functions $u_1$ and $u_2$ represent the front and rear parts of the body's surface and the functions $f_1$ and $f_2$ are determined by the temperature and composition of the medium and by the velocity of the body. A more detailed description of this model (including the study of some minimal resistance problems) can be found in \cite{temp}.
\end{zam}

Isaac Newton in his {\it Principia} \cite{N} studied the problem of minimizing the resistance in the class of rotationally symmetric bodies with fixed length and width. In modern terms, one should minimize the functional \eqref{resN} in the class of radially symmetric functions $u(x) = \phi(|x|)$, where the unknown function $\phi: [0,\, 1] \to \RRR$ is concave and monotone non-increasing, and satisfies the inequalities $0 \le \phi(\xi) \le M$. Here $M>0$ is the parameter of the problem.

The paper by Buttazzo and Kawohl \cite{BK} gave rise to studying the problem of minimal resistance in various functional classes: classes of concave (not necessarily radially symmetric) functions \cite{BFK,BrFK,BG,CLR03,K,LO,W,LZ1,JCA}, classes of functions satisfying the condition of single reflection \cite{BFK,CL1,CL2,PSIA,PSIAM,MMOP}, and some other classes \cite{BelKa,BW,LP1,K,temp,0r,Roughening,ARMA,LZ2}.

In this paper we deal with the following minimization problem. Let $\Om \subset \RRR^2$ be a convex body, that is, a compact bounded convex set with nonempty interior, and let $f : \RRR^2 \to \RRR$ be a positive continuous function. Fix a value $M > 0$. The problem is:
\beq\label{minProblem}
F(u) = \int_\Om f(\nabla u(x))\, dx \to \inf
\eeq
in the class of concave functions $u : \Om \to \RRR$ satisfying $0 \le u(x) \le M$.\footnote{Notice that since $u$ is concave, the gradient $\nabla u(x)$ exist for almost all $x$, and therefore, the integral in \eqref{minProblem} is well defined.}

The problem has been intensively studied, especially in the case when $\Om$ is a circle and $f(x) = 1/(1 + |x|^2)$. However, despite the apparent simplicity of the statement, this problem is still open. It is well known that the problem has a solution \cite{BFK,BG}.

Here we answer positively the question stated in 1993 in \cite{BK}: is it true that if $u$ is a solution to Problem \eqref{minProblem} then $u\rfloor_{\pl\Om} = 0$?

Suppose that the function $f$ satisfies the following condition:
\vspace{2mm}

 {\bf (A)}\ $\frac{(1+|x|)f(x)}{|y|f(y)} \to +\infty$ \ as\, $(1+|x|)/|y| \to 0$.
\vspace{2mm}

\hspace*{-6.5mm}Here $x = (x^1, x^2)$ and $y = (y^1, y^2)$ are points of $\RRR^2$. By $(x,z) = (x^1, x^2, z)$ we denote points of $\RRR^3$.

This condition should be understood in the usual sense: for any $N > 0$ there exists $\del > 0$ such that for all $x,\, y$ satisfying $(1+|x|)/|y| < \del$ we have $\frac{(1+|x|)f(x)}{|y|f(y)} > N$.

One easily derives that a function $f$ of the form $f(x) = c|x|^{-\al} (1 + o(1)),\, x \to \infty$ with $c>0$  satisfies Condition A, if $\al > 1$, and does not satisfy, if $0 <\al \le 1$. In particular, $f(x) = 1/(1 + |x|^2)$ satisfies this condition.

Let $x_0 \in \pl\Om$. We consider the following condition.

\vspace{2mm}

 {\bf (B)}\  $x_0$ is a regular point of the boundary, that is, the support line $l$ at $x_0$ is unique. Additionally, the endpoints of the segment $l \cap \pl\Om$ are also regular.
\vspace{2mm}

The segment $l \cap \pl\Om$ may degenerate to the point $x_0$; in this case the second part of Condition (B) is not needed.

\begin{theorem}\label{t1}
Let $u$ be a solution to Problem \eqref{minProblem}, where $f$ satisfies Condition A, and let the point $x_0 \in \pl\Om$ satisfy Condition B. Then $u(x_0) = 0$.
\end{theorem}

The following Theorem \ref{t2} is a direct consequence of Theorem \ref{t1}.

\begin{theorem}\label{t2}
Let $u$ be a solution to Problem \eqref{minProblem}, where $f$ satisfies Condition A and all points of $\pl\Om$ are regular. Then $u\rfloor_{\pl\Om} = 0$.
\end{theorem}

\section{Proof of Theorem \ref{t1}}

Assume the contrary: $u(x_0) > 0$. We are going to come to a contradiction.

Without loss of generality we assume that $u$ is upper semicontinuous. Otherwise we substitute $u$ with the function $\cl u$ defined on $\Om$ whose subgraph is the closure of the subgraph of $u$. Such a function exists and is of course unique, it is concave and upper semicontinuous,  and satisfies the inequalities $\cl u(x) \ge u(x)$ and $\cl u(x) = u(x)$ in the interior of $\Om$. As a consequence, one has $0 \le \cl u(x) \le M$ and $\cl u(x_0) \ge u(x_0) > 0$. Proofs of analogous properties for {\it convex}, rather than {\it concave}, function can be found in the book of Rockafellar \cite{R} (Part II, Section 7).

Since $u$ is upper semicontinuous and $l \cap \pl\Om$ is compact, the maximum value of $u\rfloor_{l\cap\pl\Om}$ is attained at a certain point; without loss of generality we assume that
$$
u(x_0) = \max_{x\in l\cap\pl\Om} u(x) =: z_0 > 0.
$$

Denote by $n$ the outward normal to $\pl\Om$ at $x_0$. The line $l$ is then determined by the equation $(x_0-x, n) = 0$, and $\Om$ lies in the half-plane $(x_0-x, n) \ge 0$; here and in what follows, $(\cdot\,, \cdot)$ means the scalar product.

Take $k>0$ and consider the plane of the equation $z = k(x_0-x, n)$. This plane contains the line $l \times \{0\}$, has slope $k$, and separates the domain of $u$ in the horizontal plane $\Om \times \{0\}$ and the vertical segment $\{ x_0 \} \times [0,\, z_0]$. Consider the auxiliary function $u^{(k)}(x) = \min \{ u(x),\, k(x_0-x, n) \}$; the subgraph of $u^{(k)}$ is the part of the subgraph of $u$ located below the plane.

We are going to prove that for $k$ sufficiently large, $F(u^{(k)}) < F(u)$, in contradiction with optimality of $u$.

Consider the planar convex body
$$
\Om_k = \{ x : u(x) \ge k(x_0-x, n) \}.
$$
We have $\nabla u^{(k)}(x) = -kn$ for $x$ in the interior of $\Om_k$.  Outside $\Om_k$ the function $u$ coincides with $u^{(k)}$, therefore
\beq\label{F}
F(u) - F(u^{(k)}(x)) = \int_{\Om_k} f(\nabla u(x))\, dx - \int_{\Om_k} f(u^{(k)}(x))\, dx = \int_{\Om_k} f(\nabla u(x))\, dx - f(-kn)\, |\Om_k|;
\eeq
here and in what follows, $|\cdot|$ means the area of a planar figure. It remains to show that the right hand side of this expression is positive for $k$ sufficiently large.

We are going to find a family of convex bodies $\tilde\Om_k \subset \Om_k$ satisfying the asymptotic relations
\beq\label{ur1}
\frac{1}{k}\,\sup_{x\in\tilde\Om_k} |\nabla u(x)| \to 0 \qquad \text{as} \ \ k \to +\infty;
\eeq
\beq\label{ur2}
\inf_{x\in\tilde\Om_k} \Big( \frac{|\tilde\Om_k|/|\Om_k|}{(1 + |\nabla u(x)|)/k} \Big) \ge \text{const} > 0 \quad
\text{for $k$ sufficiently large.}
\eeq
Loosely speaking, we require that first, the maximum of $|\nabla u|$ in $\tilde\Om_k$ is asymptotically much smaller than the value of $|\nabla u^{(k)}|$ in $\Om_k$ (which is equal to $k$) and second, the relative area of $\tilde\Om_k$ in $\Om_k$ decreases not too rapidly. Let us show that \eqref{ur1} and \eqref{ur2} lead to a contradiction.

The right hand side of \eqref{F} can be transformed as follows,
$$
\frac{1}{f(-kn)\, |\Om_k|} \Big( \int_{\Om_k} f(\nabla u(x))\, dx - f(-kn)\, |\Om_k| \Big) \ge
$$
$$
\frac{1}{|\tilde\Om_k|} \int_{\tilde\Om_k}
\frac{|\tilde\Om_k|/|\Om_k|}{(1 + |\nabla u(x)|)/k}\ \frac{(1 + |\nabla u(x)|) f(\nabla u(x))}{k f(-kn)}\, dx - 1 \ge
$$
\beq\label{infinf}
\inf_{x\in\tilde\Om_k} \Big( \frac{|\tilde\Om_k|/|\Om_k|}{(1 + |\nabla u(x)|)/k} \Big)
\inf_{x\in\tilde\Om_k} \Big( \frac{(1 + |\nabla u(x)|) f(\nabla u(x))}{k f(-kn)}  \Big) -1.
\eeq
According to \eqref{ur2}, the former infimum in \eqref{infinf}is greater than or equal to a positive constant. Taking for any $k>0$ a point $x_k \in \tilde\Om_k$, by \eqref{ur1} one has  $(1 + |\nabla u(x_k)|)/k \to 0$ as $k \to +\infty$, and hence, by Condition A,
$$\frac{(1 + |\nabla u(x_k)|) f(\nabla u(x_k))}{k f(-kn)} \to +\infty \qquad \text{as} \ \ k \to +\infty.$$
Taking the infima over all $x_k \in \tilde\Om_k$, one obtains that the latter infimum in \eqref{infinf} goes to infinity.
It follows that the expression in \eqref{infinf} tends to $+\infty$ as $k \to +\infty$ and therefore, the right hand side of \eqref{F} is positive for $k$ sufficiently large.

It remains to choose $\tilde\Om_k$ in such a way that \eqref{ur1} and \eqref{ur2} are satisfied. This will finish the proof of the theorem.

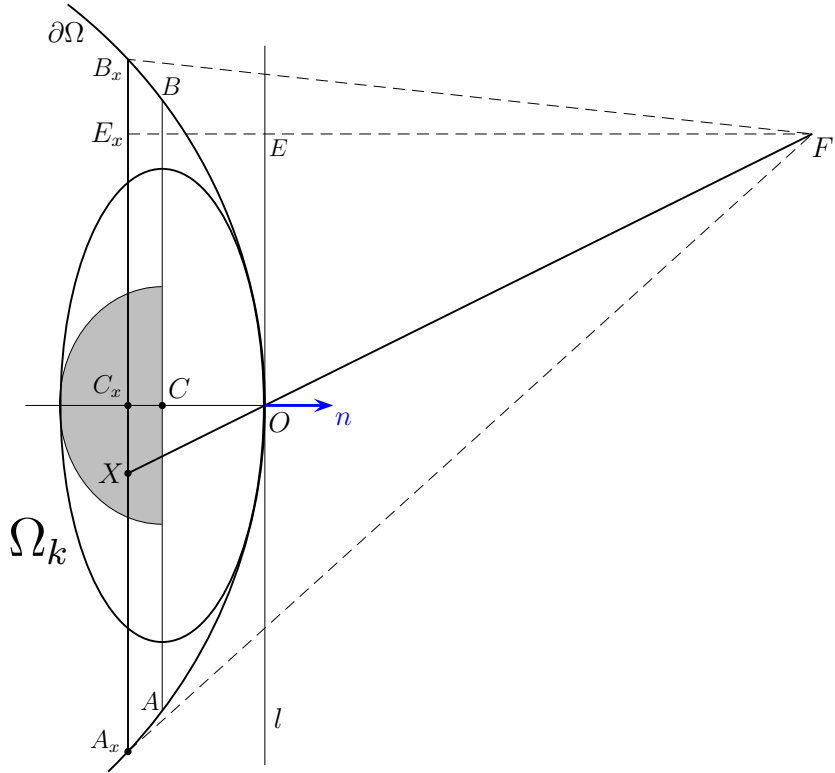
\begin{figure}[h]
\begin{picture}(0,275)
 \rput(5,4.6){
 \scalebox{0.9}{
\psellipse[fillcolor=lightgray,fillstyle=solid,linewidth=0pt](0,0)(1.5,1.75)
\pspolygon[linewidth=0pt,linecolor=white,fillcolor=white,fillstyle=solid](0,-2)(0,2)(2,2)(2,-2)

\psarc(-6,0){7.5}{-46}{52}
\psellipse(0,0)(1.5,3.5)
   \psline[linewidth=0.4pt](-2,0)(1.5,0)
\psline[linewidth=0.4pt](0,-4.5)(0,4.5)
\psline(-0.5,5.1)(-0.5,-5.1)
\psdots[dotsize=3pt](-0.5,-1)(-0.5,-5.1)(0,0)(-0.5,0)
 \psline(-0.5,-1)(9.5,4)
 \psline[linewidth=0.4pt,linestyle=dashed](-0.5,4)(9.5,4)
 \psline[linewidth=0.4pt,linestyle=dashed](-0.5,-5.1)(9.5,4)(-0.5,5.1)
 \rput(9.65,3.8){$F$}
  \rput(-0.83,-4.95){\scalebox{0.9}{$A_{x}$}}
  \rput(-0.8,4.95){\scalebox{0.9}{$B_{x}$}}
  \rput(-0.18,-4.36){\scalebox{0.9}{$A$}}
  \rput(0.13,4.7){\scalebox{0.9}{$B$}}
   \rput(-0.75,-1){$X$}
      \rput(-0.8,4){$E_{x}$}
      \rput(1.7,3.8){\scalebox{0.85}{$E$}}
      \rput(0.25,0.25){$C$}
  \rput(-0.8,0.28){\scalebox{0.9}{$C_{x}$}}
      \rput(1.72,-0.25){$O$}
      \psline[linewidth=0.4pt](1.5,5.3)(1.5,-5.3)
            \rput(1.7,-4.6){$l$}
            \psline[arrows=->,arrowscale=1.5,linewidth=1.2pt,linecolor=blue](1.5,0)(2.5,0)
            \rput(2.65,-0.2){$\blue n$}
   \rput(-1.8,-2){\scalebox{1.9}{$\Om_k$}}
      \rput(-1.4,5.5){$\pl\Om$}
  }}
  \end{picture}
\caption{In this figure, the circular arc through $A$, $B$, and $O$ is a part of $\pl\Om$, the ellipse is $\Om_k$, the shadowed domain is $\tilde\Om_k$, the points $x_0$ and $x$ are marked by $O$ and $X$, respectively. The figure corresponds to the case when $l \cap \pl\Om$ is a singleton and the point $X$ lies below the line $OC$.}
\label{fig1}
\end{figure}

\begin{figure}[h]
\begin{picture}(0,210)
 \rput(6,3.5){
 \scalebox{1.5}{

      \pspolygon[linewidth=0pt,linecolor=white,fillstyle=solid,fillcolor=lightgray]
      (1,0.24)(1.15,0.4)(1.325,0.52)(1.5,0.602)(1.5,-0.602)(1.325,-0.52)(1.15,-0.4)(1,-0.24)
\psline[linewidth=0.5pt](2,-1.5)(2,1.5)
\psarc[linewidth=0.5pt](2.5,0){1.581}{109.45}{161.58}
\psarc[linewidth=0.5pt](2.5,0){1.581}{198.42}{250.55}
\psline[linewidth=0.5pt](1,0.5)(1,-0.5)
\psline[linewidth=0.5pt](1.5,1.23)(1.5,-1.23)
     \psdots[dotsize=2.5pt](2,0)(4.1,1.2)(1.3,-0.4)
     \psecurve[linewidth=0.5pt](1.9,1.1)(2,1.5)(1.8,1.9)(0.35,2.5)(-0.9,2)(-1.6,1)(-1.6,-1)(-0.9,-2)(0.35,-2.5)(1.8,-1.9)(2,-1.5)(1.9,-1.1)
      \psline[linewidth=0.6pt,linestyle=dotted,dotsep=2pt](1.3,-1)(1.3,-2.24)(0.4,-2.24)(0.4,2.24)(1.3,2.24)(1.3,0.5)
                      \psline[linewidth=0.5pt](1.3,-0.4)(4.1,1.2)
                      \psline[linewidth=0.5pt,linestyle=dashed,dash=3pt 2pt](2,1.5)(4.1,1.2)(2,-1.5)
                      \psline[linewidth=0.5pt,linestyle=dashed,dash=3pt 2pt](1.3,-0.4)(1.3,0)(2,0)
    \rput(1.25,0.21){\scalebox{0.6}{$C_{x}$}}
    \rput(1.25,-0.65){\scalebox{0.67}{$X$}}
          \rput(1.65,0.17){\scalebox{0.67}{$C$}}
        \rput(2.18,0.27){\scalebox{0.67}{$O$}}
        \rput(4.2,1){\scalebox{0.67}{$F$}}
\pscurve[linewidth=0.5pt](1,0.24)(1.15,0.4)(1.325,0.52)(1.5,0.602)
\pscurve[linewidth=0.5pt](1,-0.24)(1.15,-0.4)(1.325,-0.52)(1.5,-0.602)
 \psline[arrows=->,arrowscale=1.5,linewidth=1pt,linecolor=blue](2,0)(2.7,0)
             \rput(2.5,-0.25){\scalebox{0.83}{$\blue n$}}
      \rput(-0.2,0){\scalebox{1.8}{$\Om$}}
      \rput(1,-1.25){\scalebox{1.1}{$\Om_k$}}
         \rput(2.1,-1.58){\scalebox{0.67}{$P$}}
         \rput(2.12,1.66){\scalebox{0.67}{$Q$}}

  }}
  \end{picture}
\caption{The case when $l \cap \pl\Om$ is a non-degenerate segment ($PQ$ in the figure) and $X$ lies below the line $OC$.}
\label{fig3}
\end{figure}
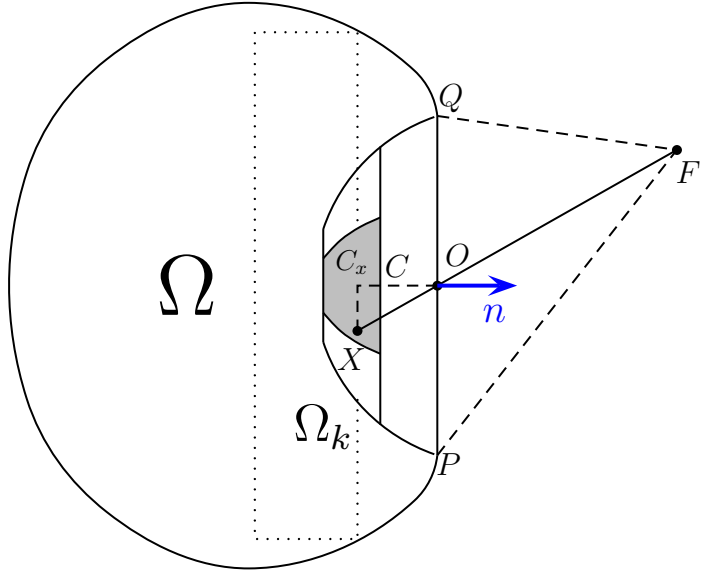

Our construction is illustrated in Figs. \ref{fig1} and \ref{fig3} corresponding to the cases when $l \cap \pl\Om$ is a point and a line segment, respectively. In both figures the point $x_0$ is indicated by the letter $O$ and the line $OC$ is orthogonal to $l$ (and its director vector is $n$).

Note that $\Om_k$ lies in the intersection of $\Om$ with the band $\{ x : 0 \le (x_0-x, n) \le M/k \}$. It follows that $\Om_k$ lies in an $\ve_k$-neighborhood of $l \cap \pl\Om$, with $\lim_{k\to+\infty} \ve_k = 0$.

Denote
$$
\al_k := \sup_{x\in\Om_k} u(x) - z_0.
$$
The sets $\Om_k,\, k >0$ form a nested family, $\Om_{k_1} \subset \Om_{k_2}$ for $k_1 \ge k_2$, therefore the function $k \mapsto \al_k$ is monotone decreasing. Since $z_0 = u(x_0)$ and $x_0 \in \Om_k$, this function is non-negative. Further, since $z_0$ is the maximal value of $u\rfloor_{l\cap\pl\Om}$ and $u$ is upper semicontinuous, one has $\al_k \to 0$ as $k \to +\infty$. It may happen than $\al_{k_0} = 0$ for a certain value $k_0$; then $\al_k$ equals zero for all $k \ge k_0$.

Consider the set
$$
\om_k := \{ (x_0-x, n) : x \in \Om_k \}.
$$
It is a closed segment contained in the positive semiaxis $[0,\, +\infty)$. Note that the orthogonal projection of $\Om_k$ on the line $x_0 + \xi n,\, \xi \in \RRR$ (the line $OC$ in Figures~\ref{fig1} and \ref{fig3}) is the line segment $x_0 - \om_k n$.

Taking $x = x_0 \in \Om_k$, one concludes that the lower endpoint of $\om_k$ is 0, and therefore, $\om_k$ has the form
$$
\om_k = [0,\, (z_0 + \bt_k)/k].
$$
Since for all $x \in \Om_k$, $(x_0-x, n) \le u(x)/k \le (z_0 + \al_k)/k$, one sees that the upper endpoint of $\om_k$ does not exceed $(z_0 + \al_k)/k$, and so, $\bt_k \le \al_k$.

Fix $\ve > 0$. For $x = x_0 - \frac{z_0 - \ve}{k} n$ one has
$$
u(x) - k(x_0-x, n) = \big[ u\big( x_0 - \frac{z_0 - \ve}{k} n \big) - z_0 \big] + \ve.
$$
For $k$ sufficiently large the expression in the right hand side of this formula is positive. This means that $x \in \Om_k$, and therefore, $(z_0 - \ve)/k$ lies in $\om_k$. Hence for $k$ sufficiently large, $\bt_k \ge -\ve$. It follows that $\bt_k \to 0$ as $k \to +\infty$.

Take $t > 0$ and consider the segment $\Om \cap \{ x : (x_0-x, n) = t \}$. It is orthogonal to the line $x_0 + \xi n,\, \xi \in \RRR$ (the line $OC$ in Figs.~\ref{fig1} and \ref{fig3}). If $t$ is sufficiently small, the segment is divided by this line into two non-degenerated segments; say the lower and the upper ones. Let the lengths of the lower and upper segments be $a(t)$ and $b(t)$. Both functions are concave, non-negative, and monotone increasing for $t$ sufficiently small.

Let a point $x$ is marked by $X$ in Fig.~\ref{fig1} and let $t = |C_{x}O| = (x_0-x,\, n)$; then $a(t) = |A_{x}C_{x}|$ and $b(t) = |B_{x}C_{x}|$.

The length of the interval $l \cap \pl\Om$ equals $a(0) + b(0)$. If the interval degenerates to the point $x_0$, one has $a(0) = b(0) = 0$. By condition B the endpoints of this interval are regular points of $\pl\Om$, hence $a'(0) = b'(0) = +\infty$. It follows that $\lim_{t\to 0} (a(t)/t) = +\infty$ and $\lim_{t\to 0} (b(t)/t) = +\infty$. We assume that $k$ is sufficiently large, so that both functions $a(t)$ and $b(t)$ are monotone increasing for $t \le (z_0 + \bt_k)/k$.

The following  formula will be needed later on. Due to concavity of $a$, for $\frac{z_0}{2k} \le \xi \le \frac{z_0 + \bt_k}{k}$ we have $a\big(\frac{z_0}{2k}\big) \ge \big( 1 - \frac{z_0}{2k\xi} \big) a(0) + \frac{z_0}{2k\xi} a(\xi)$, hence
\beq\label{estimate}
a(\xi) \le \frac{2k\xi}{z_0}\, a\Big(\frac{z_0}{2k}\Big) - \frac{2k\xi - z_0}{z_0}\, a(0) \le 2\frac{z_0 + \bt_k}{z_0}\, a\Big(\frac{z_0}{2k}\Big).
\eeq
A similar formula holds for the function $b$.

For $0 < \theta < 1$ consider the linear map $T_{\theta} : \RRR^2 \to \RRR^2$ that leaves points of the line $OC$ unchanged and moves other points in the direction orthogonal to $OC$, so as for all $x\in \RRR^2$ the distance from $T_{\theta}x$ to $OC$ is $\theta$ times the distance from $x$ to $OC$. It is defined by the equation $T_{\theta} x = \theta x + (1-\theta) ((x-x_0,n) n + x_0)$. This map is a compression with the ratio $\theta$ in the direction orthogonal to $n$.

We take a positive function $k \mapsto \theta_k$ satisfying the conditions
\beq\label{theta1}
\rm{(i)} \ \theta_k \to 0 \ \, \text{as} \ \,  k \to +\infty, \ \ \rm{(ii)} \ k\, a\big( \frac{z_0}{2k} \big)\, \theta_k \to +\infty \ \text{and} \, \ k\, b\big( \frac{z_0}{2k} \big)\, \theta_k \to +\infty \ \, \text{as} \ \,  k \to +\infty,
\eeq
and (in the case when $\al_k$ is always positive)
\beq\label{theta2}
\rm{(iii)} \ \ \frac{\theta_k}{\al_k} \to +\infty \quad \text{as} \  \, k \to +\infty.\hspace*{50mm}
\eeq
One can take, for example, $\theta_k = 1\big/\sqrt{ka\big( \frac{z_0}{2k} \big)} + 1\big/\sqrt{kb\big( \frac{z_0}{2k} \big)} + \sqrt{\al_k}$.

Let $\Om_k^+ = \Om_k \cap \{ x: (x_0-x, n) \ge\frac{ z_0}{2k} \}$, and let $\tilde\Om_k = T_{\theta_k}(\Om_k^+)$; that is, $\tilde\Om_k$ is the image of $\Om_k^+$ under the compression with the ratio $\theta_k$. In Figs. \ref{fig1} and \ref{fig3}, $\tilde\Om_k$ is shown shadowed.

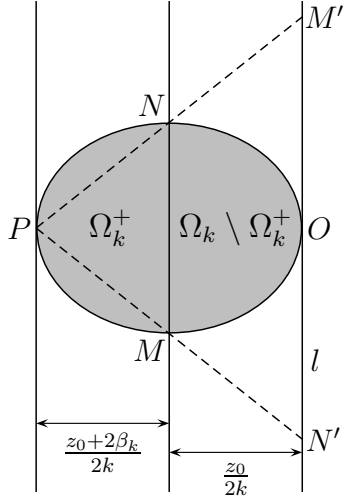
\begin{figure}[h]
\begin{picture}(0,185)
 \rput(7.5,3.6){
 \scalebox{0.7}{
\psellipse[fillcolor=lightgray,fillstyle=solid](0,0)(2.5,2)

   \psline[linewidth=0.8pt](-2.5,-5)(-2.5,4.3)
   \psline[linewidth=0.8pt](0,-5)(0,4.3)
   \psline[linewidth=0.8pt](2.5,-5)(2.5,4.3)
   \rput(-2.82,0){\scalebox{1.4}{$P$}}

 \psline[linewidth=0.8pt,linestyle=dashed](-2.5,0)(2.5,4)
 \psline[linewidth=0.8pt,linestyle=dashed](-2.5,0)(2.5,-4)
  \rput(-0.37,-2.32){\scalebox{1.4}{$M$}}
  \rput(-0.33,2.32){\scalebox{1.4}{$N$}}
    \rput(2.92,4){\scalebox{1.4}{$M'$}}
  \rput(2.9,-4){\scalebox{1.4}{$N'$}}

  \rput(-1.1,0){\scalebox{1.5}{$\Om_k^+$}}
  \rput(1.25,0){\scalebox{1.5}{$\Om_k \setminus \Om_k^+$}}

            \rput(2.75,-2.5){\scalebox{1.4}{$l$}}
            \rput(2.82,0){\scalebox{1.4}{$O$}}
            \psline[arrows=<->,arrowscale=1.7,linewidth=0.6pt](-2.5,-3.7)(0,-3.7)
            \psline[arrows=<->,arrowscale=1.7,linewidth=0.6pt](0,-4.25)(2.5,-4.25)
  \rput(-1.25,-4.25){\scalebox{1.4}{$\frac{z_0+2\bt_k}{2k}$}}
  \rput(1.25,-4.75){\scalebox{1.4}{$\frac{z_0}{2k}$}}
  }}
  \end{picture}
\caption{The sets $\Om_k$ and $\Om_k \setminus \Om_k^+$.}
\label{figAreas}
\end{figure}

The set $\Om_k^+$ is contained between the lines $(x_0-x, n) = \frac{z_0 + \bt_k}{k}$ and $(x_0-x, n) = \frac{z_0}{2k}$, and the set $\Om_k \setminus \Om_k^+$ is contained between the line $(x_0-x, n) = \frac{z_0}{2k}$ and the line $l$ of the equation $(x_0-x, n) = 0$; see Fig.~\ref{figAreas}. Let the intersection of $\Om_k$ with the line $(x_0-x, n) = \frac{z_0}{2k}$ be the segment $MN$. The intersection of $\Om_k$ with the line $(x_0-x, n) = \frac{z_0 + \bt_k}{k}$ is nonempty; take a point $P$ in this intersection. Let $M'$ and $N'$ be the points of intersection of $l$ with the lines $PM$ and $PN$, respectively.

The area of the triangle $PMN$ equals
$$
|\triangle PMN| = \frac{1}{2}\, \frac{z_0 + 2\bt_k}{2k}\, |MN|.
$$
One easily sees that $|M'N'| = \frac{2(z_0 + \bt_k)}{z_0 + 2\bt_k}\, |MN|$, and therefore, the area of the trapezoid $MNN'M'$ is
$$
|\Box MNN'M'| = \frac{z_0}{2k}\, \Big( \frac{1}{2} + \frac{z_0 + \bt_k}{z_0 + 2\bt_k} \Big) |MN|.
$$

The set $\Om_k^+$ contains the triangle $PMN$, and the set $\Om_k \setminus \Om_k^+$ is contained in the trapezoid $MNN'M'$, therefore
$$
\frac{|\Om_k \setminus \Om_k^+|}{|\Om_k^+|} = \frac{|\Om_k|}{|\Om_k^+|} - 1 \le \frac{z_0(3z_0 + 4\bt_k)}{(z_0 + 2\bt_k)^2} \quad \Rightarrow \quad
\frac{|\Om_k|}{|\Om_k^+|} \le \frac{4(z_0 + \bt_k)^2}{(z_0 + 2\bt_k)^2}.
$$

Thus, one has
\beq\label{Om}
\frac{|\tilde\Om_k|}{|\Om_k|} = \theta_k \frac{|\Om_k^+|}{|\Om_k|} \ge \frac{\theta_k}{4}\, \frac{(z_0 + 2\bt_k)^2}{(z_0 + \bt_k)^2}
= \frac{\theta_k}{4}\, (1 + o(1)) \quad \text{as} \ \, k \to +\infty.
\eeq

Let $x \in \tilde\Om_k$ be a regular point of $u$; it is indicated by letter $X$ in Figs.~\ref{fig1}  and \ref{fig3}. If the tangent plane to the graph of $u$ at $x$ is not horizontal, then its intersection with the horizontal plane $z=0$ is a straight line, say $l_x$, and
\beq\label{nab}
|\nabla u(x)| = \frac{u(x)}{\text{dist}(x,l_x)}.
\eeq
Of course $l_x$ does not intersect the interior of $\Om$.

\begin{figure}[h]
\begin{picture}(0,120)
 \rput(-1.2,0.7){
 \scalebox{1}{
\pspolygon[linewidth=0pt,linecolor=white,fillstyle=solid,fillcolor=lightgray](1,3.37)(1.5,3.2)(2,3)(3.5,2.18)(4.5,1.57)(4.75,1.35)(4.98,1.1)(5,1)(5,0)(1,0)
\psecurve[linewidth=0.5pt,linewidth=1.2pt,linecolor=brown](0.5,3.5)(1,3.37)(1.5,3.2)(2,3)(3.5,2.18)(4.5,1.57)(4.98,1.1)(5,1)(4.97,0.9)
 \psline[arrows=->,arrowscale=1.5,linewidth=1pt](2,0)(9,0)
\psline(1,3.5)(8.5,-0.25)
\psline[linestyle=dashed,linewidth=0.5pt](2,0)(2,3)
\psline(5,0)(5,1)
\psdots(2,3)(5,1)(2,0)(6.4,0)
\pspolygon[linewidth=0pt,linecolor=white,fillstyle=solid,fillcolor=white](1.6,-0.1)(1.6,3.75)(0,3.75)(0,-0.1)
\rput(2,-0.27){$X$}
\rput(5,-0.25){$O$}
\rput(6.4,-0.28){$F$}
\rput(2.75,3.3){$(x,u(x))$}
\rput(5.66,0.75){$(x_0,z_0)$}
\rput(5,1.8){$\sigma_x$}
\rput(1.2,-0.65){(i)}
 }}
 \rput(9.5,0.7){
 \scalebox{1}{
 \pspolygon[linewidth=0pt,linecolor=white,fillstyle=solid,fillcolor=lightgray](0,0.9)(1,1.5)(1.75,1.9)(2.5,2.15)(3,2.3)(3.5,2.3)(3.8,2.24)(4,2)(4,0)(0,0)
 \psecurve[linewidth=0.5pt,linewidth=1.2pt,linecolor=brown](-0.5,0.5)(0,0.9)(1,1.5)(2.5,2.15)(3.5,2.3)(4,2)(3.7,1.75)
 \psline(0,1)(5,3.5)
 \rput(1,-0.27){$X$}
\rput(4,-0.25){$O$}
\psdots(1,1.5)(4,2)(1,0)
 \psline[linestyle=dashed,linewidth=0.5pt](1,0)(1,1.5)
 \psline(4,0)(4,2)
 \psline[arrows=->,arrowscale=1.5,linewidth=1.2pt](1,0)(5.6,0)
 \pspolygon[linewidth=0pt,linecolor=white,fillstyle=solid,fillcolor=white](0.7,-0.1)(0.7,1.4)(-0.5,1.4)(-0.5,-0.1)
 \rput(0.2,1.8){$(x,u(x))$}
\rput(4.7,2){$(x_0,z_0)$}
\rput(3.5,3.05){$\sigma_x$}
 \rput(0,-0.65){(ii)}
 }}
  \end{picture}
\caption{The section of graph$(u)$ and the tangent plane to graph$(u)$ at $x$ by the vertical plane through $X$ and $O$. (i) The line $\sigma_x$ intersects the ray $XO$; (ii) $\sigma_x$ does not intersect this ray.}
\label{fig05}
\end{figure}
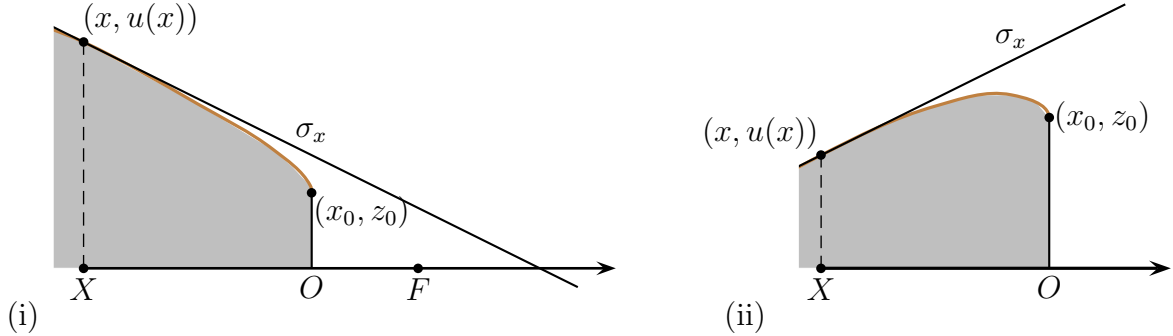

The intersection of this tangent plane with the vertical plane through the line $XO$ (and therefore, through the points $(x,0)$ and $(x_0,0)$) is a straight line, say $\sigma_x$, that contains the point $(x,u(x))$ and lies above the point $(x_0,u(x_0))$; see Fig.~\ref{fig05}. The slope of this line in the direction $\overrightarrow{XO}$ is greater than or equal to $(u(x_0) - u(x))/|x_0-x|$.

It may happen that (i) the line $\sigma_x$ intersects the ray with the vertex at $x$ and with the director vector $x_0-x$ (the ray $XO$ in Figs.~\ref{fig1}, \ref{fig3}, and \ref{fig05}), or (ii) it does not intersect this ray.
If the case (i) is realized, one has $u(x) > u(x_0)$ (and therefore $\al_k > 0$), and the slope of the line $\sigma_x$ is greater than or equal to $-\al_k/|x_0-x|$. It follows that the distance between $x_0$ and the point of intersection of $\sigma_x$ with the ray $XO$ is greater than or equal to $\frac{z_0}{\al_k} |x_0-x|$. As a consequence, the point of intersection lies on the ray behind the point $x_0 + \frac{z_0}{\al_k} (x_0-x)$ (the point $F$ in Figs.~\ref{fig1}, \ref{fig3} and \ref{fig05}). For the point $F$ we have the proportion
\beq\label{proportion}
|OF| = \frac{z_0}{\al_k}\, |XO|.
\eeq
Hence $l_x$ does not intersect the interior of the convex hull of $\Om$ and the segment $[X F]$,
\beq\label{intersect}
l_x \cap \text{int}\big(\text{Conv}(\Om \cup [XF])\big) = \emptyset.
\eeq

If the case (ii) is realized then $l_x$ does not intersect the interior of the convex hull of the union of $\Om$  and the ray $XO$. In particular, formula \eqref{intersect} remains true, where in the case $\al_k = 0$ the point $F$ can be imagined as the infinitely remote point on the ray and $[XF]$ should be understood as the ray $XO$.

Let $A_{x}$ and $B_{x}$ be the endpoints of the segment $\{ x' : (x', n) = (x, n) \} \cap \Om$ (see Fig.~\ref{fig1}), and let $\pl_x \Om$ be the part of $\pl\Om$ on the left of the line $A_{x}B_{x}$. In other words, $\pl_x\Om := \pl\Om \cap \{ x' : (x', n) \le (x, n) \}$.  Let $[A_{x}F]$ and $[B_{x}F]$ denote the corresponding segments, $A_{x}F$ and $B_{x}F$ be the lines containing these segments, and $[A_{x}FB_{x}]$ be the union of the segments. (If $\al_k = 0$, $[A_{x}F]$ and $[B_{x}F]$ should be understood as the rays co-directional with the ray $XO$ with the vertices at $A_{x}$ and $B_{x}$, respectively, and $[A_{x}FB_{x}]$ as the union of these rays.)

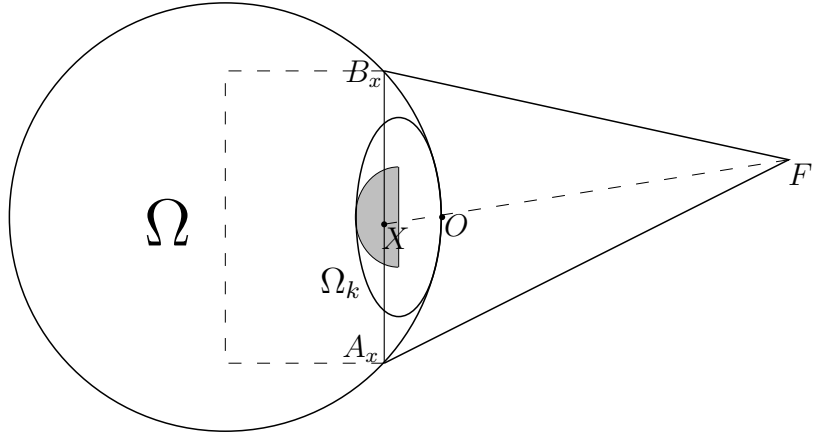
\begin{figure}[h]
\begin{picture}(0,150)
 \rput(8.5,2.5){
 \scalebox{0.38}{
\psellipse[fillcolor=lightgray,fillstyle=solid,linewidth=0pt](0,0)(1.5,1.75)
\pspolygon[linewidth=0pt,linecolor=white,fillcolor=white,fillstyle=solid](0,-2)(0,2)(2,2)(2,-2)
\psline[linewidth=0.1pt](0,-1.75)(0,1.75)

\pscircle[linewidth=1.5pt](-6,0){7.5}
\psellipse[linewidth=1.5pt](0,0)(1.5,3.5)
\psline[linewidth=1.1pt](-0.5,5.1)(-0.5,-5.1)
\psdots[dotsize=6pt](-0.5,-0.25)(1.5,0)
 \psline[linestyle=dashed,dash=12.5pt 15pt,linewidth=1pt](-0.5,-0.25)(13.5,2)
 \psline[linewidth=1.5pt](-0.5,-5.1)(13.5,2)(-0.5,5.1)
 \rput(13.9,1.5){\scalebox{2.6}{$F$}}
  \rput(-1.3,-4.6){\scalebox{2.6}{$A_{x}$}}
  \rput(-1.2,5){\scalebox{2.6}{$B_{x}$}}
   \rput(-0.1,-0.75){\scalebox{2.6}{$X$}}
      \rput(2,-0.3){\scalebox{2.6}{$O$}}
   \rput(-2,-2.3){\scalebox{3}{$\Om_k$}}
      \rput(-8,-0.2){\scalebox{5.5}{$\Om$}}
\psline[linestyle=dashed,dash=12.5pt 15pt](-0.5,5.1)(-6,5.1)(-6,-5.1)(-0.5,-5.1)


  }}
  \end{picture}
\caption{The domain bounded by the curve $\pl_x\Om$ and the broken line $[A_{x}FB_{x}]$. In this figure, $\pl_x\Om$ is the arc of circumference with the endpoints $A_{x}$ and $B_{x}$ located to the left of the line $A_{x}B_{x}$.}
\label{fig2}
\end{figure}

The open domain bounded by the union of the curve $\pl_x \Om$ and the broken line (or the union of rays) $[A_{x}FB_{x}]$ (see Fig.~\ref{fig2}) is contained in Conv$(\Om \cup [XF])$, and therefore, does not intersect the line $l_x$. Therefore we have
$$
\text{dist}(x, l_x) \ge \text{dist}(x, \pl_x\Om \cup [A_{x}FB_{x}]) = \min \{ \text{dist}(x, \pl_x\Om),\, \text{dist}(x, [A_{x}F]),\, \text{dist}(x, [B_{x}F]) \}.
$$

Take $t_0$ sufficiently small, so as the functions $a(t)$ and $b(t)$ are monotone increasing for $0 \le t \le t_0$. Consider the rectangle located on the left of the line $A_{x}B_{x}$ such that the segment $[A_{x}B_{x}]$ is one of its vertical sides (see Figs.~\ref{fig2} and \ref{fig3}) and the length of its horizontal sides equals $t_0/2$. One easily sees that for $k$ sufficiently large (namely, for $k$ satisfying $(z_0 + \bt_k)/k \le t_0/2$) and for $x \in \tilde\Om_k$ this rectangle is contained in $\Om$.

The distance from $x$ to the union of the three sides of the rectangle (except $A_{x}B_{x}$) equals $\min \{ t_0/2,\, |XA_{x}|,\, |XB_{x}| \}$. It follows that
$$
\text{dist}(x, \pl_x\Om) \ge \min \{ t_0/2,\, |xA_{x}|,\, |xB_{x}| \} \ge \min \{ t_0/2,\, \text{dist}(x, [A_{x}F]),\, \text{dist}(x, [B_{x}F]) \}.
$$
Hence
\beq\label{dist}
\text{dist}(x, l_x) \ge \min \{ t_0/2,\, \text{dist}(x, [A_{x}F]),\, \text{dist}(x, [B_{x}F]) \}.
\eeq

Let $C = C^k$ indicate the point $x_0 - \frac{z_0}{2k}n$.    
Denote by $C_{x}$ the projection of $X$ on the line $CO$, and by $E_{x}$ and $E$ the projections of $F$ on the lines $A_{x}B_{x}$ and $l$, respectively\footnote{If $\al_k = 0$, the points $E_{x}$ and $E$ are not defined.} (see Fig.~\ref{fig1}). Note that $|CO| = \frac{z_0}{2k}$.

Assuming that $\al_k > 0$, using \eqref{proportion} and taking into account that the triangles $XC_{x}O$ and $OEF$ are similar, one obtains
\beq\label{eq1}
|E_{x}F| \ge |EF| = |C_{x}O|\, \frac{|OF|}{|XO|} \ge \frac{z_0}{2k}\, \frac{z_0}{\al_k} = \frac{z_0^2}{2\al_k k}.
\eeq

Denote $\xi = |C_{x}O|$; we have $\frac{z_0}{2k} \le \xi \le \frac{z_0 + \bt_k}{k}$ and
\beq\label{ac2}
|A_{x}C_{x}| = a(\xi) \ge a\big( \frac{z_0}{2k} \big).
\eeq
Using formula \eqref{estimate}, one obtains
\beq\label{ac}
|A_{x}C_{x}| \le 2\, \frac{z_0 + \bt_k}{z_0}\, a\big( \frac{z_0}{2k} \big).
\eeq

Consider two cases.
\vspace{2mm}

{\bf (a)}  $X$ lies below the line $OC$, that is, in the same half-plane bounded by $OC$ as $A_{x}$.

First assume that (i) the line $\sigma_x$ intersects the ray $XO$. Using \eqref{ac} we obtain
$$
|XC_{x}| \le \theta_k |A_{x}C_{x}| \le 2\theta_k\, \frac{z_0 + \bt_k}{z_0}\, a\big( \frac{z_0}{2k} \big).
$$
It follows that
\beq\label{eq2}
\frac{|XE_{x}|}{|E_{x}F|} = \frac{|XC_{x}|}{|C_{x}O|} \le \frac{|XC_{x}|}{|CO|} \le 4\theta_k\, \frac{ka\big( \frac{z_0}{2k} \big)}{z_0}\ \frac{z_0 + \bt_k}{z_0}.
\eeq

Further, using \eqref{ac2} one has
\beq\label{eq3}
|A_{x} X| = |A_{x} C_{x}| - |C_{x} X| \ge (1-\theta_k) |A_{x}C_{x}| \ge (1 - \theta_k) a\big( \frac{z_0}{2k} \big).
\eeq

Let us evaluate the distance $h$ between $X$ and the line $A_{x} F$. Consider the triangle $A_{x} X F$. The area $S$ of this triangle can be calculated in two ways,
$$
S = \frac{1}{2}\, |A_{x} X|\, |E_{x}F| \quad \text{and} \quad S = \frac{1}{2}\, h\, |A_{x} F| \le \frac{1}{2}\, h\, (|A_{x} X| + |XE_{x}| + |E_{x}F|),
$$
and using inequalities \eqref{eq1}, \eqref{eq2}, and \eqref{eq3} one obtains
$$
h \ge \frac{|A_{x} X|\, |E_{x}F|}{|A_{x} X| + |XE_{x}| + |E_{x}F|} = \frac{1}{\frac{1}{|E_{x}F|} + \frac{1}{|A_{x} X|} \big( 1 + \frac{|XE_{x}|}{|E_{x}F|} \big)}
$$
$$
\ge \frac{1}{\frac{2\al_k k}{z_0^2} + \frac{1}{1 - \theta_k}\ \frac{1}{a\big( \frac{z_0}{2k} \big)} \Big( 1 + 4\theta_k\, \frac{ka\big( \frac{z_0}{2k} \big)}{z_0}\, \frac{z_0 + \bt_k}{z_0} \Big)}
 = \frac{1}{k\theta_k}\ \frac{1}{\frac{2\al_k}{z_0^2 \theta_k} + \frac{1}{1 - \theta_k}\ \frac{1}{ka\big( \frac{z_0}{2k} \big) \theta_k} + \frac{4}{1 - \theta_k}\, \frac{z_0 + \bt_k}{z_0^2}}.
$$
It follows that
\beq\label{eq4}
\text{dist}(x, [A_{x}F]) \ge \text{dist}(x, A_{x}F) \ge \frac{1}{k\theta_k}\, H_k^1,
\eeq
where
\beq\label{eq5}
H_k^1 = \frac{1}{\frac{2\al_k}{z_0^2\theta_k} + \frac{1}{1 - \theta_k}\Big( \frac{1}{ka\big( \frac{z_0}{2k} \big)\theta_k} + \frac{1}{kb\big( \frac{z_0}{2k} \big)\theta_k} \Big) + \frac{4}{1 - \theta_k}\, \frac{z_0 + \bt_k}{z_0^2}}.
\eeq
Note that $H_k^1$ is chosen to be invariant with respect to exchanging $a$ and $b$. Due to \eqref{theta1} and \eqref{theta2}, the denominator in this expression tends to $4/z_0$, hence $H_k^1 \to z_0/4$ as $k \to +\infty$.

Now assume that (ii) $\sigma_x$ does not intersect the ray $XO$.  If $\al_k > 0$ then the argument above in the case (a) remains valid, and formulae \eqref{eq4} and \eqref{eq5} hold true. If $\al_k = 0$, the argument is valid when $F$ is replaced with an arbitrary point $x_0 + \frac{z_0}{\al}\, (x_0 - x)$, $\al > 0$ on the ray $XO$. It follows that inequality \eqref{eq4} is true when $\al_k$ is replaced with arbitrary $\al > 0$ in \eqref{eq5}. Taking the limit $\al \to 0$, one concludes that \eqref{eq4} is also true when we have $\al = \al_k = 0$ in \eqref{eq5}.
\vspace{2mm}

{\bf (b)} $X$ lies above the line $OC$. Here again we consider two cases.

\begin{figure}[h]
\begin{picture}(0,170)
 \rput(2,3.25){
 \scalebox{0.5}{
\psellipse[fillcolor=lightgray,fillstyle=solid,linewidth=0pt](0,0)(1.5,1.75)
\pspolygon[linewidth=0pt,linecolor=white,fillcolor=white,fillstyle=solid](0,-2)(0,2)(2,2)(2,-2)

\psarc(-6,0){7.5}{-46}{52}
\psellipse(0,0)(1.5,3.5)
   \psline[linewidth=0.8pt](-2,0)(1.5,0)
\psline[linewidth=1.6pt](0,-4.5)(0,4.5)
\psline[linewidth=1.6pt](-0.5,5.1)(-0.5,-6.6)
\psdots[dotsize=5pt](-0.5,1.5)(-0.5,-5.1)(0,0)(-0.5,0)(1.5,0)(9.5,-6)
 \psline(-0.5,1.5)(9.5,-6)
       \psline[linewidth=1.6pt,linestyle=dotted,dotsep=6pt](-0.5,-6)(9.5,-6)
 \psline[linewidth=0.8pt,linestyle=dashed,dash=10pt 6pt](-0.5,-5.1)(9.5,-6)(-0.5,5.1)
 \rput(10,-5.6){\scalebox{2}{$F$}}
  \rput(-1.1,-4.8){\scalebox{1.8}{$A_{x}$}}
  \rput(-1,4.8){\scalebox{1.8}{$B_{x}$}}
   \rput(-0.9,1.9){\scalebox{2}{$X$}}
      \rput(-1,-6.1){\scalebox{1.8}{$E_{x}$}}
      \rput(0.45,-0.45){\scalebox{2}{$C$}}
  \rput(-1,-0.4){\scalebox{1.8}{$C_{x}$}}
      \rput(1.9,0.35){\scalebox{2}{$O$}}
      \rput(-2.2,5.6){\scalebox{2}{$\pl\Om$}}
            \rput(-3.4,-5.6){\scalebox{2}{(b$_1$)}}
  }}
  \rput(11,3.25){
 \scalebox{0.5}{
\psellipse[fillcolor=lightgray,fillstyle=solid,linewidth=0pt](0,0)(1.5,1.75)
\pspolygon[linewidth=0pt,linecolor=white,fillcolor=white,fillstyle=solid](0,-2)(0,2)(2,2)(2,-2)

\psarc(-6,0){7.5}{-46}{52}
\psellipse(0,0)(1.5,3.25)
   \psline[linewidth=0.8pt](-2,0)(1.5,0)
\psline[linewidth=1.6pt](0,-4.5)(0,4.5)
\psline[linewidth=1.6pt](-0.5,5.1)(-0.5,-5.1)
\psdots[dotsize=5pt](-0.5,0.6)(-0.5,-5.1)(0,0)(-0.5,0)(1.5,0)(9.5,-2.4)(-0.5,-2.4)
 \psline(-0.5,0.6)(9.5,-2.4)
 \psline[linewidth=1.6pt,linestyle=dotted,dotsep=6pt](-0.5,-2.4)(9.5,-2.4)
 \psline[linewidth=0.8pt,linestyle=dashed,dash=10pt 6pt](-0.5,-5.1)(9.5,-2.4)(-0.5,5.1)
 \rput(9.9,-2){\scalebox{2}{$F$}}
  \rput(-1.1,-4.8){\scalebox{1.8}{$A_{x}$}}
  \rput(-1,4.8){\scalebox{1.8}{$B_{x}$}}
   \rput(-0.9,1){\scalebox{2}{$X$}}
      \rput(-1.2,-2.5){\scalebox{1.8}{$E_{x}$}}
      \rput(0.45,-0.45){\scalebox{2}{$C$}}
  \rput(-1,-0.4){\scalebox{1.8}{$C_{x}$}}
      \rput(1.9,0.35){\scalebox{2}{$O$}}
      \rput(-2.2,5.6){\scalebox{2}{$\pl\Om$}}
                  \rput(-3.4,-5.6){\scalebox{2}{(b$_2$)}}
  }}
  \end{picture}
\caption{$X$ lies above the line $OC$. The two figures correspond to the cases (b$_1$) when $E_{x}$ lies below $A_{x}$ and (b$_2$) when $E_{x}$ lies between $C_{x}$ and $A_{x}$.}
\label{fig4}
\end{figure}
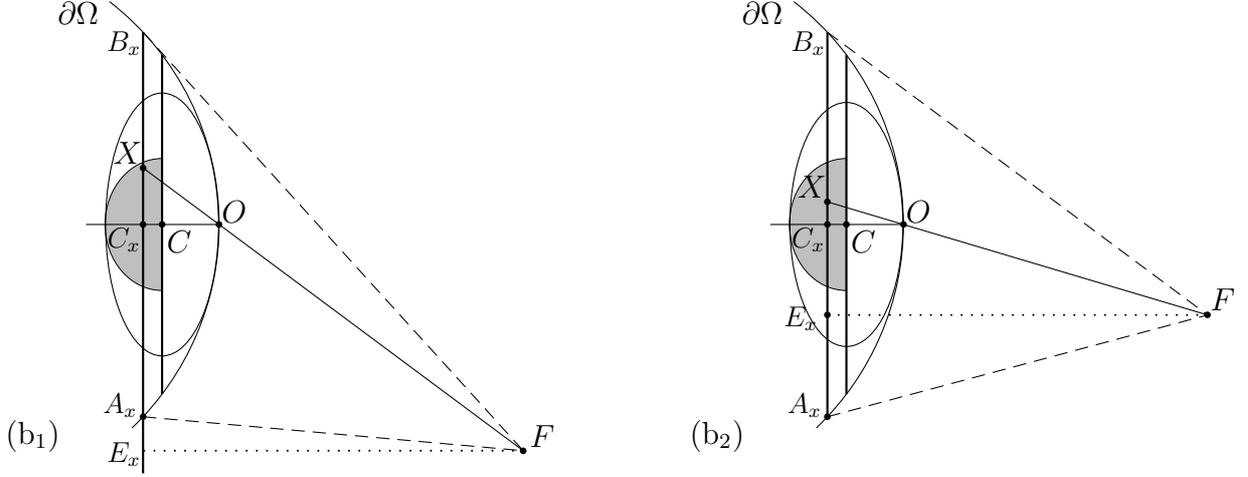

(b$_1$) $E_{x}$ lies below $A_{x}$ or (in the case $\al_k = 0$) does not exist; see Fig.~\ref{fig4}\,(b$_1$). In this case, using \eqref{ac2} and \eqref{eq5}, we have
$$
\text{dist}(x, [A_{x}F]) = |XA_{x}| \ge |A_{x}C_{x}| = a(\xi) \ge a\big( \frac{z_0}{2k} \big) > \frac{1}{k\theta_k}\, H_k^1.
$$

(b$_2$) $E_{x}$ lies between the points $A_{x}$ and $C_{x}$; see Fig.~\ref{fig4}\,(b$_2$). Due to \eqref{ac} and \eqref{ac2}, one has
\beq\label{ax}
|A_{x}E_{x}| \le |A_{x}C_{x}| \le 2\, \frac{z_0 + \bt_k}{z_0}\, a\big( \frac{z_0}{2k} \big) \qquad
\text{and} \qquad
|A_{x}X| \ge |A_{x}C_{x}| \ge a\big( \frac{z_0}{2k} \big).
\eeq
Let us consider again the area $S$ of the triangle $A_{x} X F$. One has
$$
S = \frac{1}{2}\, |A_{x} X|\, |E_{x}F| \quad \text{and} \quad S = \frac{1}{2}\, h\, |A_{x} F| \le \frac{1}{2}\, h\, (|A_{x} E_{x}| + |E_{x}F|),
$$
and using \eqref{eq1} and \eqref{ax} one obtains
$$
h \ge \frac{|A_{x} X|\, |E_{x}F|}{|A_{x} E_{x}| + |E_{x}F|} = \frac{|A_{x} X|}{\frac{|A_{x} E_{x}|}{|E_{x}F|} + 1}
\ge \frac{a\big( \frac{z_0}{2k} \big)}{2\, \frac{z_0 + \bt_k}{z_0}\, a\big( \frac{z_0}{2k} \big) \frac{2\al_k k}{z_0^2} + 1} \ge \frac{1}{k\theta_k}\, H_k^2,
$$
where
$$
 H_k^2 =  \frac{1}{\frac{4\al_k}{\theta_k}\, \frac{z_0 + \bt_k}{z_0^3} + \frac{1}{ka\big( \frac{z_0}{2k} \big) \theta_k} + \frac{1}{kb\big( \frac{z_0}{2k} \big) \theta_k}}  \to +\infty \ \ \text{as}  \  k \to +\infty.
$$
$H_k^2$ is also invariant with respect to exchanging $a$ and $b$.

Thus, in all cases one has
$ \text{dist}(x, [A_{x}F]) \ge \frac{1}{k\theta_k} H_k,$
where $H_k = \min\{ H_k^1,\, H_k^2 \} \to z_0/4$ as $k \to +\infty$.

In a completely similar way one derives the estimate $ \text{dist}(x, [B_{x}F]) \ge \frac{1}{k\theta_k} H_k$. Using that the function $a$ is bounded and that by \eqref{theta1} $k\, a\big( \frac{z_0}{2k} \big) \theta_k$ goes to infinity, we have
$$
\frac{1}{k\theta_k} = \frac{a\big( \frac{z_0}{2k} \big)}{k\, a\big( \frac{z_0}{2k} \big) \theta_k} \to 0 \quad \text{as}  \  k \to +\infty,
$$
hence by \eqref{dist} for $k$ large enough we have
$$
\text{dist}(x, l_x) \ge \min \Big\{ \frac{t_0}{2},\ \frac{1}{k\theta_k} H_k \Big\} = \frac{1}{k\theta_k} H_k.
$$
By \eqref{nab},  for a regular $x \in \tilde\Om_k$ one has
\beq\label{nabla}
\frac{1}{k}\, |\nabla u(x)| = \frac{1}{k}\ \frac{u(x)}{\text{dist}(x,l_x)} \le \frac{M}{H_k}\, \theta_k  \to 0 \ \ \text{as}  \  k \to +\infty.
\eeq
Thus, equation \eqref{ur1} is true. Further, using \eqref{Om} and \eqref{nabla} one obtains
$$
\frac{|\tilde\Om_k|/|\Om_k|}{(1 + |\nabla u(x)|)/k} \ge \frac{1}{4}\, \frac{(z_0 + 2\bt_k)^2}{(z_0 + \bt_k)^2}\  \frac{1}{\frac{1}{k\theta_k} + \frac{M}{H_k}}
= \frac{z_0}{16M}\, (1 + o(1)) \quad \text{as}  \  k \to +\infty.
$$
Thus, equation \eqref{ur2} is also true. This completes the proof of the theorem.

\section*{Acknowledgements}

{This work was supported by Foundation for Science and Technology (FCT), within project UID/MAT/04106/2019 (CIDMA). I am very grateful to G. Wachsmuth for a discussion, which helped to fix an error in the text.

\end{document}